\documentclass[letterpaper, 10 pt, conference]{ieeeconf}  

\usepackage{subfigure}
\usepackage{graphics} 
\usepackage{epsfig} 
\usepackage{mathptmx} 
\usepackage{times} 
\usepackage{amsmath} 
\usepackage{amssymb}  
\usepackage{makeidx}
\usepackage{float}
\usepackage{balance}
\usepackage{booktabs}
\usepackage[ruled]{algorithm2e}
\usepackage{algorithmic}
\usepackage{bbm}
\usepackage{bm}

\IEEEoverridecommandlockouts
\title{\LARGE \bf
Constrained Sampling-based Trajectory Optimization using Stochastic Approximation
}

\author{George I. Boutselis, Ziyi Wang and Evangelos A. Theodorou
\thanks{The authors are with the school of Aerospace Engineering, Georgia Institute of Technology, USA
        {\tt\small gbouts, ziyiwang, evangelos.theodorou@gatech.edu}.}%
}

\begin{document}

\maketitle
\thispagestyle{empty}
\pagestyle{empty}

\begin{abstract}

We propose a sampling-based trajectory optimization methodology for constrained problems. We extend recent works on stochastic search to deal with box control constraints, as well as nonlinear state constraints for discrete dynamical systems. Regarding the former, our strategy is to optimize over truncated parameterized distributions on control inputs. Furthermore, we show how non-smooth penalty functions can be incorporated into our framework to handle state constraints. Simulations on cartpole and quadcopter show that our approach outperforms previous methods on constrained sampling-based optimization, in terms of quality of solutions and convergence speed.

\end{abstract}

\section{INTRODUCTION}
Numerical trajectory optimization algorithms have received increasing attention over the past few years, due to their applicability to a wide range of problems in reinforcement learning, controls and robotics. Generally speaking, existing approaches can be classified into two main categories: gradient- and sampling-based methods. Gradient-based methods use local approximations of the cost and dynamics functions to approach a solution and have demonstrated successful results in simulated and real scenarios \cite{tassa2014control}. Popular implementations rely on direct (constrained) optimization techniques \cite{kobilarov2008discrete} or dynamic programming principles \cite{li2004iterative}.

Gradient information allows the above methods to attain fast convergence. Nonetheless, in many practical applications gradient information may not be available, and thus a different optimization strategy needs to be considered. These include problems with inherent discontinuities, such as contact dynamics \cite{tassa2012synthesis}, or environments where only evaluations of the costs and dynamics are provided \cite{todorov2012mujoco}. Potentially, one can use smoothing techniques and/or learn data-based models to provide gradient information, however numerical instabilities and poor generalizations may often arise \cite{pereira2019learning}.

Due to the aforementioned limitations, sampling-based techniques have been extensively used over the past few years. Even though convergence is typically slower than their gradient-based counterparts, their ability to handle discontinuities and lack of gradient information is higly desirable. Popular implementations include path integral control \cite{theodorou2010generalized} and the cross entropy method \cite{rubinstein2013cross}, which have recently found application on real robotics settings \cite{buchli2011learning}. Loosely speaking, these schemes compute control updates based on the evaluation of the costs associated with sampled trajectories.

Despite their popularity, extensions of sampling-based control methods to constrained problems are limited. The simplest and most common approach to handle state constraints is to associate infeasible trajectories with high costs, usually in the form of indicator functions \cite{kalakrishnan2011stomp}. In this way, control updates tend to iteratively avoid infeasible regions of the state space. Box control constraints are typically handled via the same approach, or even by applying clamping on the controls sampled \cite{tassa2014control}. Recently, a constrained version of cross entropy was suggested in \cite{wen2018constrained} for reinforcement learning problems, which utilized an extra elitist sorting of the rollouts with respect to their associated magnitude of constraint violation.

In this paper, we develop a constrained, sampling-based trajectory optimizer by relying on the theory of stochastic approximation (SA) \cite{robbins1951stochastic}, \cite{spall2005introduction}. SA methods tackle the problem by optimizing over parameterized distributions of the decision variables, and subsequently sampling candidate solutions from them \cite{hu2011stochastic}, \cite{zhou2014gradient}. This approach is applied iteratively until convergence criteria have been met.

Our contributions lie in: (i) formulating discrete trajectory optimization as a SA problem, (ii) handling box control constraints directly by considering truncated parameterized distributions and (iii) incorporating non-smooth penalty functions to account for state constraints. Extension (ii) arises from our choice to optimize over parametric distributions, and thus no further steps are required to deal with box constraints. Moreover, non-smooth penalty terms have been shown to possess important advantages over their smoothed counterparts \cite{nocedal2006numerical}, and, as shown in this work, they can naturally be used within our sampling-based framework. Numerical simulations on robotics tasks and different optimization scenarios show that our approach outperforms previous works on constrained sampling control.

The remainder of this paper is organized as follows: in Section \ref{section:unconstrained} we introduce the stochastic optimization scheme in an unconstrained setting. We extend the framework to problems with box control constraints and nonlinear state constraints in Section \ref{section:constrained}. Our proposed algorithm is compared against other sampling-based approaches in Section \ref{section:simulation}. Finally, we conclude the paper and discuss future research directions in Section \ref{section:conclusions}.

\section{Unconstrained stochastic optimization}
\label{section:unconstrained}

\subsection{Derivation for discrete dynamic systems}
Our sampling control method has been inspired by the work in \cite{zhou2014gradient}. Therein, the authors attempt to minimize costs without structural properties such as differentiability and convexity. Here, we modify the particular framework to account for discrete-time Markovian dynamics.

Let us begin with the stochastic, discrete-time optimal control problem:
\begin{equation}
\label{original_problem}
\begin{split}
\min_{\bm{\theta}}\hspace{1mm} &\mathbb{E}\big[J\big(\bm{x}, \bm{u}(\bm{x};\bm{\theta})\big)\big]\\
\text{s.t.}\hspace{1.8mm} x^{k+1}=&f^k(x^k, u^k(x^k;\theta^k)),\hspace{1.8mm} k=0, 1, ..., H.
\end{split}
\end{equation}
Here, $x^k\in\mathbb{R}^n$, $u^k\in\mathbb{R}^m$, $\theta ^k\in\mathbb{R}^d$ denote the state, control input and associated control parameters, respectively, at the $k$th time instance. That is, the control inputs $u^k$ will be parameterized by a set of parameter vectors $\theta^k$, over which we will be optimizing the expected cost $\mathbb{E}[J]$. The stochastic transition dynamics $f^k$ impose a probability density function (p.d.f.) for the next states, $x^{k+1}\sim p(x^{k+1}|x^k, u^k)$. We will denote the corresponding state, control and parameter sequences as $\bm{x}:=(({x^{0}})^\top, ..., ({x^H})^\top)$,  $\bm{u}:=(({u^{0}})^\top, ..., ({u^{H-1}})^\top)$ and $\bm{\theta}:=(({\theta^{0}})^\top, ..., ({\theta^{H-1}})^\top)$, respectively.

\begin{table}
\centering
  \begin{tabular}{| l | c | }
    \hline
    $x^k$, $\bm{x}$ & States \\
    $u^k$, $\bm{u}$ & Control inputs \\
    $\theta^k$, $\bm{\theta}$ & Control parameters \\
    $\chi^k(\cdot)$, $\chi(\cdot)$ & Probability density function of $\bm{\theta}$
    \\
    $\rho^k$, $\bm{\rho}$ & Parameters of p.d.f. $\chi(\bm{\theta})$\\
    \hline
  \end{tabular}
  \caption{Variable notation}
\end{table}

Inspired by \cite{zhou2014gradient}, we will introduce a sampling distribution for the control parameters $\bm{\theta}$, over which we will optimize the expected value of $J(\cdot)$. Let us denote by $\chi(\bm{\theta};\bm{\rho})$ this distribution of $\bm{\theta}$ which will be parameterized by the vector $\bm{\rho}:=(({\rho^{0}})^\top, ..., ({\rho^{H-1}})^\top)$. Then, we will aim to minimize:

\begin{equation}
\label{min_exp}
\begin{split}
\min_{\bm{\rho}} \hspace{1mm}&\mathbb{E}\big[J\big(\bm{x}, \bm{u}(\bm{x};\bm{\theta})\big)\big],\\
	\text{s.t.}\quad &x^{k+1}=f(x^k, u^k(\theta^k)),\\
		&\theta^{k}\sim \chi^k(\theta^k; \rho^{k}),\quad k=0, 1, ..., H,
\end{split}
\end{equation}
We will henceforth assume that $\chi(\cdot)$ belongs to the exponential family of distributions \cite{brown1986fundamentals}.

To justify this approach, notice that the optimum cost of \eqref{min_exp} is always an upper bound to the minimum cost in \eqref{original_problem}, associated with $\bm\theta_*$. Moreover, the two costs become equal when the entire probability mass of $\chi(\cdot)$ is concentrated on $\bm{\theta}_*$. Hence, it is intuitive to minimize the expected value of the original cost with respect to a distribution over $\bm{\theta}$.

Further, notice that we can allow $J$ to be non-convex and even discontinuous, since we optimize with respect to $\bm{\rho}$ which only appears in $\chi(\cdot)$. Our controller, $u^k(\theta^k)$, also provides great flexibility since it can define a feedforward, feedback, or even non-linear policy. These advantages stem from our choice to optimize \eqref{min_exp} not directly with respect to the control variables $\bm\theta$, but in terms of the parameters $\bm\rho$ of their sampling distributions. A sketch of the corresponding optimization methodology is described by Algorithm \ref{alg:methodology}.

\begin{algorithm} \label{alg:methodology}
 \KwData{Initial parameters $\bm{\rho}$}
 \While{$\bm{\rho}$ has not converged}{
   Sample candidate control parameters $\bm{\theta}$ from $\chi(\bm{\theta};\bm{\rho})$ and evaluate their costs $J$\;
   Use a gradient-based method to update $\bm{\rho}$ and thus $\chi(\bm{\theta};\bm{\rho})$\;
 }
 \caption{Methodology for stochastic optimization}
\end{algorithm}

{\it Example:} One simple example of this formulation is the following: Suppose we parameterize our control with a linear policy $u^k=K^kx+c^k$, where now $\theta^k$ corresponds to $\{K^k, c^k\}$. If we consider each $K^k$ and $c^k$ to follow Gaussian distributions, then \eqref{min_exp} is minimized over their corresponding means and variances.

To proceed, we will introduce a continuous {\it shape} function, $S:\mathbb{R}\rightarrow\mathbb{R}^{+}$, such that $S(y)$ is monotonically decreasing in $y$, $\forall y\in\mathbb{R}$. One example is $S(y)=\exp(-\kappa y)$ for some $\kappa>0$, which we will also use in our simulations. We will also take a logarithmic transformation of the expected cost in \eqref{min_exp} before performing minimization. We introduce these auxiliary transforms because (i) it has been shown that certain selections empirically improve numerical implementation \cite{zhou2014gradient}, and (ii) is is easier to show connections with Path Integral control-related schemes. Therefore, we will attempt to solve the equivalent problem:
\begin{equation}
	\begin{split}
		\label{stoch_opt_prob_new}
		\max_{\bm{\rho}}&\mathcal{L}(\bm{\rho})=\max_{\bm{\rho}}\ln\big( \mathbb{E}\big[S\big(J\big(\bm{x}, \bm{u}(\bm{\theta})\big)\big)\big]\big)\\
		=&\max_{\bm{\rho}}\ln\bigg( \int S\big(J(\bm{\theta})\big)\Gamma(\bm\theta)\prod_{k=0}^{H} \chi^k(\theta^{k}; \rho^{k})\mathrm{d}\theta^{0}\cdots\mathrm{d}\theta^{H}\bigg).
	\end{split}
\end{equation}
where we have dropped the implicit dependence of $x^k$ and $u^k$ on $\theta^k$. The last line of \eqref{stoch_opt_prob_new} is obtained by assuming that the transition dynamics $f^k(\cdot)$ are Markovian, as well as imposing independence of $\theta^k$'s with respect to other terms. Hence, the distribution of a state/control sequence can be written as:
\begin{equation*}
	\begin{split}
		&p(\bm{x}, \bm{u},\bm{\theta})=p(x^{0})p(u^{0}|x^0;\theta^{0})\chi^0(\theta^0;\rho^0)p(x^1|u^0, x^0)	\cdots\\ &p(u^{H-1}|x^{H-1};\theta^{{H-1}})\chi^{H-1}(\theta^{H-1};\rho^{H-1})p(x^H|u^{H-1}, x^{H-1}))\\
		&\hspace{14.7mm}=\Gamma(\bm\theta)\prod_{k=0}^{H} \chi^k(\theta^{k}; \rho^{k}),
	\end{split}
\end{equation*}
where $\Gamma(\cdot)$ is defined accordingly. Now, to maximize equation \eqref{stoch_opt_prob_new}, we will have to compute the derivatives of $\mathcal{L}$ with respect to $\bm{\rho}$. Since each $\chi^k$ belongs to the exponential family, we can write \cite{brown1986fundamentals}
\begin{equation}
    \label{eq:expon}
    \chi^k(\theta^{k};\rho^k)=\exp((\rho^{k})^{\top}T(\theta^{k})-\phi(\rho^{k})),
\end{equation}were $T(\theta^{k})$ denotes the vector of sufficient statistics and $\phi(\rho^{k}):=\ln(\int\exp((\rho^{k})^{\top}T(\theta^{k}))\mathrm d\theta^{k})$. By pushing the gradient inside the integral in \eqref{stoch_opt_prob_new}, one can explicitly compute \cite{zhou2014gradient}
\begin{equation}
	\label{grad}
	\nabla_{\rho^k}\mathcal{L}=\mathbb{E}_{\mathcal{P}}[T(\theta^{k})] - \mathbb{E}_{\chi^{k}}[T(\theta^{k})],
\end{equation}
where $\mathbb{E}_{\mathcal{P}}$ denotes expectation with respect to the ``path" distribution defined by
\begin{equation}
\label{pathprob}
\mathcal{P}:=\frac{S\big(J(\{\theta^k\})\Gamma(\{\theta^k\})\prod_{k=0}^{H} \chi^k(\theta^{k}; \rho^{k})}{\int S\big(J(\{\theta^k\})\Gamma(\{\theta^k\})\prod_{k=0}^{H} \chi^k(\theta^{k}; \rho^{k})\mathrm{d}\theta^{0}\cdots\mathrm{d}\theta^{H}}
\end{equation}
and $\mathbb{E}_{\chi^{k}}$ denotes expectation with respect to the distribution of parameters $\theta^{k}$. Similarly, one can define higher-order derivatives. Based on the expressions above, a gradient ascent scheme for \eqref{stoch_opt_prob_new} simply reads:
\begin{equation}
	\label{update_grad}
	(\rho^k)_{i+1}=(\rho^k)_{i} + \gamma_{i}\big(\mathbb{E}_{\mathcal{P}_{i}}[T(\theta^{k})] - \mathbb{E}_{\chi^{k}_i}[T(\theta^{k})]\big),
\end{equation}
where $i$ corresponds to the $i$th iteration of the algorithm and $\gamma_{i}$ is a (possibly) iteration-dependent learning rate. The leftmost expectation above is computed through sampling, while the remaining term can be computed analytically. More complex updates can be designed in an analogous manner. By changing the natural parameters $\rho$ in an iterative fashion, we expect to converge to a distribution of the control variables $\theta$ that minimize our cost.

\subsection{Comparison with previous works}
The developed optimization algorithm can be viewed as a generalization of stochastic control schemes derived from Path Integral control theory \cite{theodorou2010generalized}. Therein, the dynamics are treated as a stochastic differential equation where Wiener noise usually enters through the control channel. To obtain a numerical algorithm, such works approximate the dynamics with an Euler-Maruyama scheme and proceed by minimizing the Kullback-Leibler of a cost between uncontrolled and controlled dynamics. Specifically, the transition dynamics are restricted to the form: $f(x^k, \bar{u}^k) = F(x^k)+G(x^k)(\bar{u}^k+\xi^k)$, where $F$, $G$ are properly defined matrices, $\xi^k\sim\mathcal{N}(0,\Delta t)$ and $\Delta t$ is the time step. The update scheme from \cite{theodorou2015nonlinear} reads:
\begin{equation}
	\label{pi}
	(\bar{u}^{k})_{i+1}=(\bar{u}^{k})_{i}+ \mathbb{E}_{\mathcal{Q}_i}[\xi^{k}],
\end{equation}
with $\mathcal{Q}:=\frac{\exp(-\kappa J(\{\theta^k\})\prod_{k=0}^{H-1} p(x^{k+1}|u^k, \theta^{k})}{\int \exp(-\kappa J(\{\theta^k\}))\prod_{k=0}^{H-1} p(x^{k+1}|u^k, \theta^{k})\mathrm{d}x^{0}\cdots\mathrm{d}x^{H}}$, $\kappa>0$.

The above expressions can essentially be viewed as a specific case of \eqref{update_grad}. Indeed, this is true when $S(y)\equiv \exp(-\kappa y)$, $\theta^k\equiv u^k$, $\gamma=\sqrt{\sigma}$ and each $\theta^k$ is Gaussian with fixed variance. To see the latter, notice that the last term in \eqref{grad} is a constant number with respect to the p.d.f. $\mathcal{Q}$. Hence, due to the Gaussianity of $\theta^k$ one has
\[\mathbb{E}_{\mathcal{Q}_{i}}[T(\theta^{k})] - \mathbb{E}_{\chi^{k}_i}[T(\theta^{k})]=\mathbb{E}_{\mathcal{Q}_{i}}[\sqrt{\sigma}\xi].\]
Similar expressions are given in works \cite{theodorou2010generalized}, where the policy has a specific form.

We stress again though that expressions \eqref{grad}, \eqref{update_grad} hold for the statistics of {\it any type of parameterized policy}. In contrast, optimizing over generic policies by following the approach in \cite{theodorou2015nonlinear} is not straightforward. Specifically, one will have to differentiate the change of measure between uncontrolled and controlled dynamics. Regularly, such computations cannot be carried out analytically, and are only specific to each candidate policy.

Similar update equations are also used within the Cross Entropy method \cite{de2005tutorial}, where only a set of elite sampled policies are used, whose number has to be pre-specified. We omit further details due to space limitations, and will make numerical comparisons with these methods in section \ref{section:simulation}.

\section{Constrained sampling-based trajectory optimization}
\label{section:constrained}

\subsection{Sampling control with box constraints}
\label{section:box_constraint}
We will now show an extension of the methodology above that accounts for box constraints on the control parameters
\[\bm l\leq\bm\theta\leq\bm u\] in a direct fashion; that is, without casting them as soft constraints.

We recall that our approach to solving \eqref{original_problem} was to iteratively sample candidate solutions from a distribution $\chi(\bm{\theta};\bm{\rho})$, and then update the underlying distribution. The idea here is to constrain ourselves to truncated distributions over the specified region. In this way, all candidate solutions will be directly sampled within the acceptable domain.

In this direction, we will let the sampling distribution of the control parameters $\bm{\theta}$ have the form:
\begin{equation}
\label{trunc_dist}
\chi(\bm{\theta};\bm{\rho})=\frac{\tilde{\chi}(\bm{\theta};\bm{\rho})}{\int_{\bm{l}}^{\bm{u}}\tilde{\chi}(\bm{\theta};\bm{\rho})\mathrm d\bm{\theta}},
\end{equation}
where $\tilde{\chi}(\bm{\theta};\bm{\rho})$ is of exponential type and $(\bm{l}, \bm{u}]$ denotes the hard control bounds\footnote{In case the lower bound can be attained, we subtract a small positive number from $\bm l$.}. Equation \eqref{trunc_dist} implies that $\chi$ is the truncated distribution corresponding to $\tilde\chi$. We will consider
\begin{equation*}
\tilde{\chi}(\bm{\theta};\bm{\rho})=\begin{cases}
\exp(\bm{\rho}^\top T(\theta) - \phi(\bm{\rho})) &\text{,  $\bm{\theta}\in (\bm{l}, \bm{u}]$}\\
0 &\text{, otherwise}
\end{cases}
\end{equation*}
where $\phi(\cdot)$ is a properly defined function \cite{brown1986fundamentals}. It is easy to verify that \eqref{trunc_dist} defines a proper distribution. Now denote for brevity the expected cost from \eqref{stoch_opt_prob_new} as $L:=\mathbb{E}[S(J)]$. We will compute the gradient of $\ln L(\bm{\rho})$ when the expectation is computed with respect to the truncated distribution \eqref{trunc_dist}. Specifically, we will have
\begin{equation*}
\nabla_{\bm{\rho}}\ln L(\bm{\rho})=\resizebox{.8\hsize}{!}{$\frac{\int_{\bm{l}}^{\bm{u}}S(J(\bm{x}, \bm{u}(\bm{x};\bm{\theta})))\Gamma(\bm{\theta})\nabla_{\bm{\rho}}\ln({\chi}(\bm{\theta};\bm{\rho})){\chi}(\bm{\theta};\bm{\rho})\mathrm d\bm{\theta}}{\int_{\bm{l}}^{\bm{u}}S(J(\bm{x}, \bm{u}(\bm{x};\bm{\theta})))\Gamma(\bm\theta){\chi}(\bm{\theta};\bm{\rho})\mathrm d\bm{\theta}}.$}
\end{equation*}
From \eqref{trunc_dist} one obtains
\begin{equation*}
\nabla_{\bm{\rho}}\ln({\chi}(\bm{\theta};\bm{\rho}))=\underbrace{\nabla_{\bm{\rho}}\big(\ln(\tilde{\chi}(\bm{\theta};\bm{\rho}))\big)}_{\Gamma_1} - \underbrace{\nabla_{\bm{\rho}}\big(\ln\big(\textstyle\int_{\bm{l}}^{\bm{u}}\tilde{\chi}(\bm{\theta};\bm{\rho})\mathrm d\bm{\theta}\big)\big)}_{\Gamma_2},
\end{equation*}
where from the previous derivation we have
\begin{equation*}
\Gamma_1=T(\bm \theta) - \mathbb{E}_{\tilde{\chi}}[T(\bm \theta)],
\end{equation*}
while
\begin{equation*}
\begin{split}
&\Gamma_2=\frac{\int_{\bm{l}}^{\bm{u}}\nabla_{\bm{\rho}}\tilde{\chi}(\bm{\theta};\bm{\rho})\mathrm d\bm{\theta}}{\int_{\bm{l}}^{\bm{u}}\tilde{\chi}(\bm{\theta};\bm{\rho})\mathrm d\bm{\theta}}=\frac{\int_{\bm{l}}^{\bm{u}}\nabla_{\bm{\rho}}\ln(\tilde{\chi}(\bm{\theta};\bm{\rho}))\tilde{\chi}(\bm{\theta};\bm{\rho})\mathrm d\bm{\theta}}{\int_{\bm{l}}^{\bm{u}}\tilde{\chi}(\bm{\theta};\bm{\rho})\mathrm d\bm{\theta}}=\\
&\frac{\int_{\bm{l}}^{\bm{u}}(T(\theta) - \mathbb{E}_{\tilde\chi}[T(\theta)])\tilde{\chi}(\bm{\theta};\bm{\rho})\mathrm d\bm{\theta}}{\int_{\bm{l}}^{\bm{u}}\tilde{\chi}(\bm{\theta};\bm{\rho})\mathrm d\bm{\theta}}=\mathbb{E}_{\chi}[T(\theta)] - \mathbb{E}_{\tilde{\chi}}[T(\theta)].
\end{split}
\end{equation*}
Hence
\begin{equation*}
\nabla_{\bm{\rho}}\ln({\chi}(\bm{\theta};\bm{\rho}))=T(\bm \theta) - \mathbb{E}_{\chi}[T(\bm \theta)],
\end{equation*}
and
\begin{equation}
\label{updatetrunc}
\nabla_{\bm{\rho}}\ln L(\bm{\rho})=\mathbb{E}_{\mathcal{P}}[T(\bm{\theta})] - \mathbb{E}_{\chi}[T(\bm{\theta})].
\end{equation}
where $\mathbb{E}_{\mathcal{P}}[T(\bm{\theta})]$ is computed numerically based on \eqref{pathprob} and $\mathbb{E}_{\chi}[T(\bm{\theta})]$ can be computed analytically for certain p.d.f.'s, including the truncated normal distribution.

\subsection{Sampling control with nonlinear state constraints}
\label{section:state_constraint}
Suppose now we have the following generic, constrained optimization problem
\begin{equation}
\label{fully_constrained_problem}
\begin{split}
\min_{\bm{\theta}}\hspace{1mm} \mathbb{E}\big[J&\big(\bm{x}, \bm{u}(\bm{x};\bm{\theta})\big)\big]\\
\text{s.t.}\hspace{1.8mm} \mathbb{E}\big[g^k_i(\bm{x}&, \bm{u}(\bm{x};\bm{\theta})\big]\leq 0,\hspace{1.8mm}i=1, ..., L\\
\mathbb{E}\big[h^k_j(\bm{x}&, \bm{u}(\bm{x};\bm{\theta})\big]= 0, \hspace{1.8mm} j=1, ..., D\\
l^k\leq&\theta^k\leq u^k,\hspace{1.8mm}k=0, ..., H,
\end{split}
\end{equation}
where we have omitted the dynamics constraints for simplicity. A standard approach for solving \eqref{fully_constrained_problem} with sampling control is to sample trajectories as in section \ref{section:unconstrained} and simply assign high costs to infeasible rollouts \cite{kalakrishnan2011stomp}. In contrast, we will explore here {\it non-smooth penalty functions}.

First, let us use the methodology of sections \ref{section:unconstrained} and \ref{section:box_constraint} and transform \eqref{fully_constrained_problem} into
\begin{equation}
\label{fully_constrained_problem_rho}
\begin{split}
\min_{\bm{\rho}}\hspace{1mm} \mathbb{E}\big[J&\big(\bm{x}, \bm{u}(\bm{x};\bm{\theta})\big)\big]\\
\text{s.t.}\hspace{1.8mm} \mathbb{E}\big[g^k_i(\bm{x}&, \bm{u}(\bm{x};\bm{\theta})\big]\leq 0,\hspace{1.8mm}i=1, ..., L\\
\mathbb{E}\big[h^k_j(\bm{x}&, \bm{u}(\bm{x};\bm{\theta})\big]= 0, \hspace{1.8mm} j=1, ..., D\\
\bm\theta\sim&\chi(\bm\theta;\bm\rho),
\end{split}
\end{equation}
where $\chi$ is the truncated distribution from \eqref{trunc_dist}. Notice that problem \eqref{fully_constrained_problem_rho} can be viewed as a deterministic optimization problem with respect to $\bm\rho$. Then, we can use well-established results from optimization theory and rewrite \eqref{fully_constrained_problem_rho} as \cite{nocedal2006numerical}
\begin{equation}
\label{expectations_barrier_rho}
\begin{split}
\min_{\bm{\rho}}\hspace{1mm} &\bigg(\mathbb{E} \big[J\big(\bm{x}, \bm{u}(\bm{x};\bm{\theta})\big)\big]+\zeta\sum_{i,k}\big(\mathbb{E} \big[g_i^k(\bm{x}, \bm{u}(\bm{x};\bm{\theta}))\big]\big)^+\\
&\hspace{13mm}+\zeta\sum_{j,k}\big|\mathbb{E} \big[h_j^k(\bm{x}, \bm{u}(\bm{x};\bm{\theta})\big]\big|\bigg)\\
&\text{s.t.}\hspace{1.8mm} 
\bm\theta\sim\chi(\bm\theta;\bm\rho),
\end{split}
\end{equation}
where $(f)^+:=\max(0, f)$ for any function $f$ and $\zeta>0$ is an external parameter which controls how the severity of constraint violation is penalized. As shown in \cite{nocedal2006numerical}, a (local) solution to \eqref{fully_constrained_problem_rho} is also a (local) minimizer to \eqref{expectations_barrier_rho} under mild assumptions, when the external parameter satisfies $\zeta>\zeta_*$, for some $\zeta_*>0$. This is because the corresponding barrier functions are {\it exact}. This contrasts with schemes where non-exact penalty functions are used (such as the quadratic barrier function \cite{nocedal2006numerical}) where $\zeta$ above has to approach infinity in order to obtain an equivalent solution to \eqref{fully_constrained_problem_rho}.

One difficutly with solving \eqref{expectations_barrier_rho} is that the utilized penalty terms are non-differentiable at zero. To get around this issue, we observe that $(\mathbb{E}[f])^+\leq\mathbb{E}[(f)^+]$ and $|\mathbb{E}[f]|\leq\mathbb{E}[|f|]$ for any (integrable) function $f$. Based on this inequality, we will thus be minimizing the upper bound of \eqref{expectations_barrier_rho}:
\begin{equation}
\label{expectations_barrier_rho_new}
\begin{split}
\min_{\bm{\rho}}\hspace{1mm} &\bigg(\mathbb{E} \bigg[J\big(\bm{x}, \bm{u}(\bm{x};\bm{\theta})\big)+\zeta\sum_{i,k}\big(g_i^k(\bm{x}, \bm{u}(\bm{x};\bm{\theta}))\big)^+\\
&\hspace{13mm}+\zeta\sum_{j,k}\big|h_j^k(\bm{x}, \bm{u}(\bm{x};\bm{\theta})\big|\bigg]\bigg)\\
&\text{s.t.}\hspace{1.8mm} 
\bm\theta\sim\chi(\bm\theta;\bm\rho).
\end{split}
\end{equation}
Now, problem \eqref{expectations_barrier_rho_new} can be solved as suggested by sections \ref{section:unconstrained} and \ref{section:box_constraint}. In particular, we can apply the shape function $S(\cdot)$ and logarithmic transform from eq. \eqref{stoch_opt_prob_new} and use the update formulas \eqref{update_grad}, \eqref{updatetrunc}.


\subsection{Algorithm}
\label{section:algorithm}

\begin{algorithm}[t]
\caption{Constrained sampling-based controller}
 \KwData{Initial parameters $\bm{\rho}_0$, initial penalty coefficient $\zeta$, constraint violation threshold $\xi$, system dynamics $f$, inequality constraint $g$, equality constraint $h$, cost function $J$, sample size $N$, step size $\alpha$, penalty increase ratio $\beta$;}
 Initialize $x^0$;\\
 Set $i=1$;\\
 \While{constraints not satisfied}{
 \While{$\bm{\rho}$ not converged}{
 \For{n = 1:N}{
 Sample control parameter trajectories $\bm{\theta}_n \sim \bm{\chi}(\bm{\theta};\bm{\rho}_{i-1})$;
 \\
 \For{k = 1:H}{
 Compute state trajectories $x_n^{k}=f(x_n^{k-1},u_n^{k-1}(x_n^{k-1};\theta_n^{k-1}))$;
 }
 Calculate penalized cost $\tilde{J}_n=J(\bm{x}_n,\bm{u}_n(\bm{x}_n;\bm{\theta}_n))+\zeta\sum_{k=1}^{H}(g(\bm{x}^k_n,\bm{u}^k_n(\bm{x}^k_n;\bm{\theta}^k_n))^+ +|h(\bm{x}^k_n,\bm{u}^k_n(\bm{x}^k_n;\bm{\theta}^k_n))|)$;
 }
 Update parameters $\bm{\rho}^{(m)}=\bm{\rho}^{(m-1)} + \alpha\Big[\sum_{n=1}^N\frac{S(\tilde{J}_n)(T(\bm{\theta}_n)-\frac{1}{N}\sum_{n=1}^NT(\bm{\theta}_n))}{\sum_{n=1}^N S(\tilde{J}_n)}\Big]$;\\
 \If{any constraint exceeds $\xi$}{
   Break\;
   }
 }
 Increase penalty coefficient $\zeta = \beta\zeta$;\\
 }
         
\label{alg:state_constraint}
\end{algorithm}

We now introduce our constrained sampling-based controller, as summarized in Algorithm \ref{alg:state_constraint}. Given the constrained optimization problem of the form \eqref{fully_constrained_problem}, we first initialize the parameters $\bm{\rho}^{(0)}$ and penalty coefficient $\zeta$. At every iteration $m$, $N$ control parameter trajectories are sampled from the constrained distribution ${\chi}(\cdot)$. The state trajectories are obtained from propagating the discrete-time dynamics forward. The constraint-penalized cost for each trajectory is calculated by adding the corresponding penalty terms to the original cost function. The parameters are updated using the Monte-Carlo approximation of a gradient ascent scheme.

As commonly done in constraint optimization, we sequentially increase the penalty coefficient $\zeta$ until we converge to a set of parameters $\bm\rho$ that satisfy the constraints. Note that having a high penalty term from the beginning usually makes the problem highly nonlinear and poses difficulties for optimizers, especially when the solution is near the constraint boundaries \cite{nocedal2006numerical}. Similarly, numerical issues arise when we let the optimizer approach highly infeasible regions for a small $\zeta$, and thus a stopping criterion has to be imposed followed by an increase in the value of $\zeta$.

\begin{figure}
\centering
\includegraphics[width=175 pt]{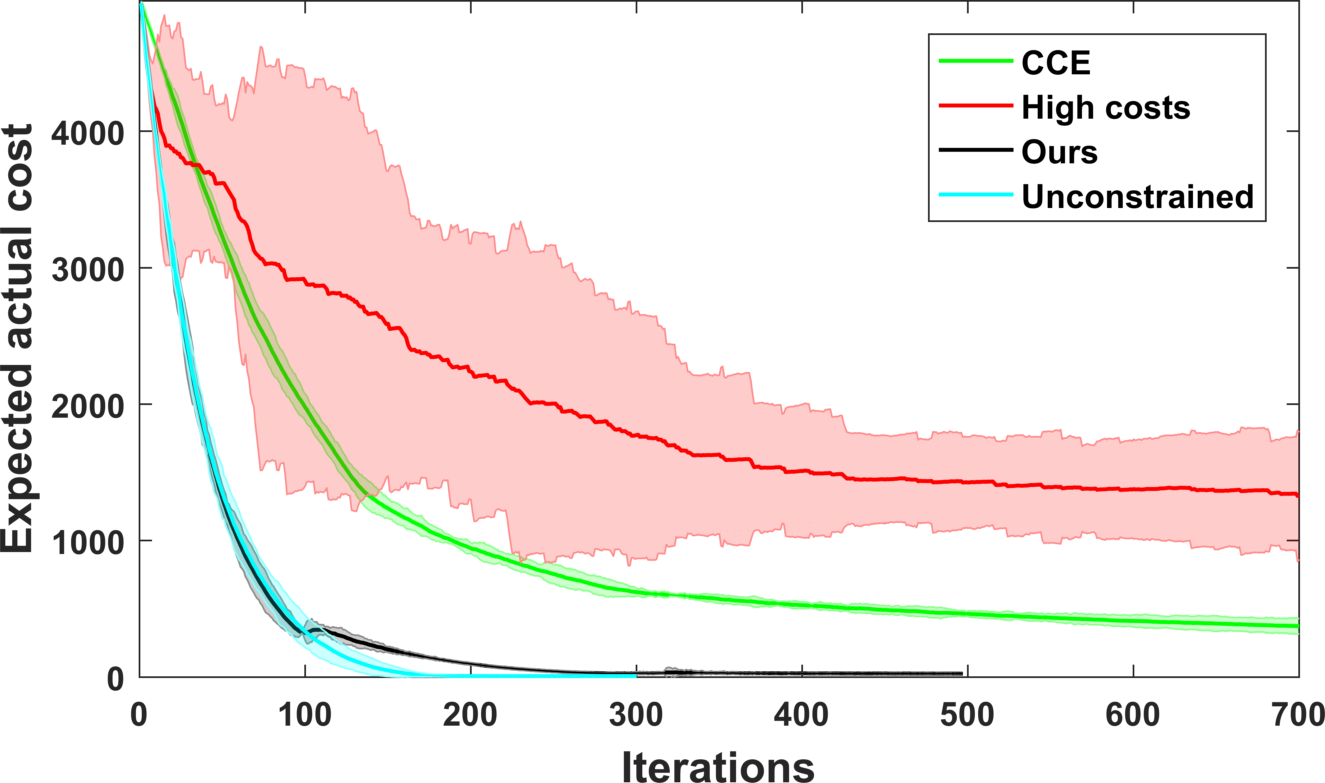}
\caption{Expected cost for cartpole with a $\pm3\sigma$ confidence interval.}
\end{figure}

\begin{figure}
\centering
\includegraphics[width=175 pt]{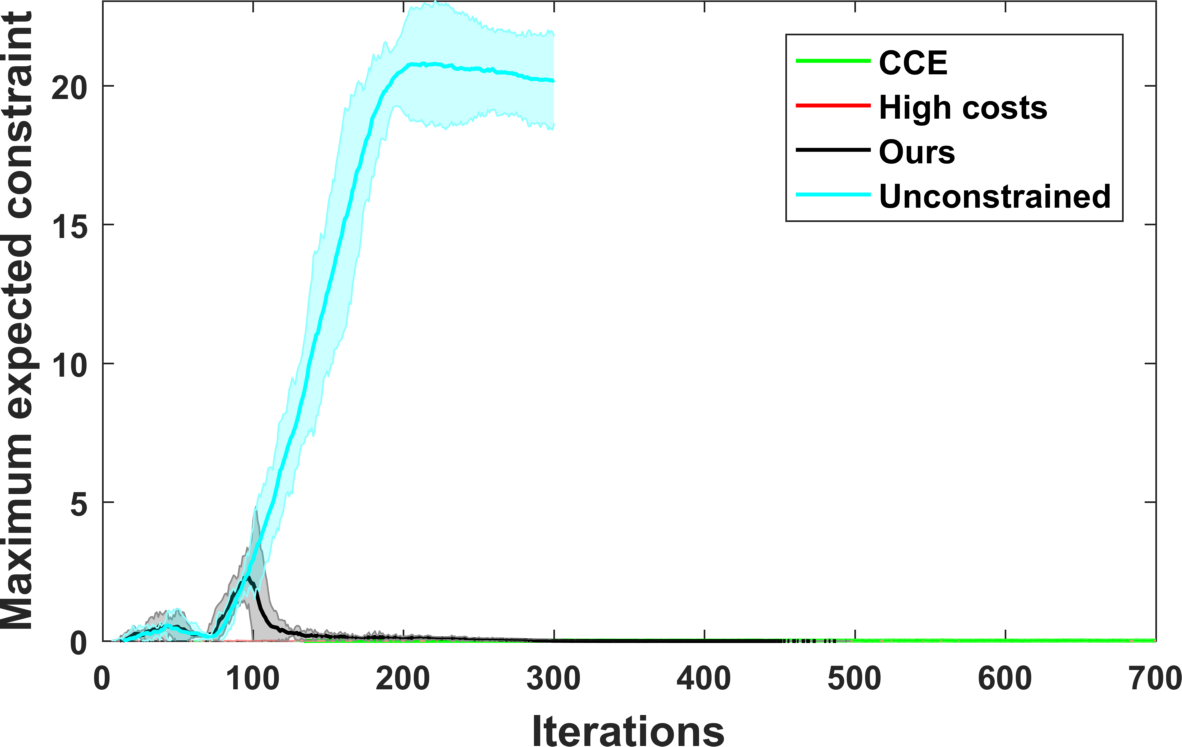}
\caption{Maximum expected constraints for cartpole with a $\pm3\sigma$ confidence interval.}
\end{figure}

\begin{figure}
\centering
\includegraphics[width=175 pt]{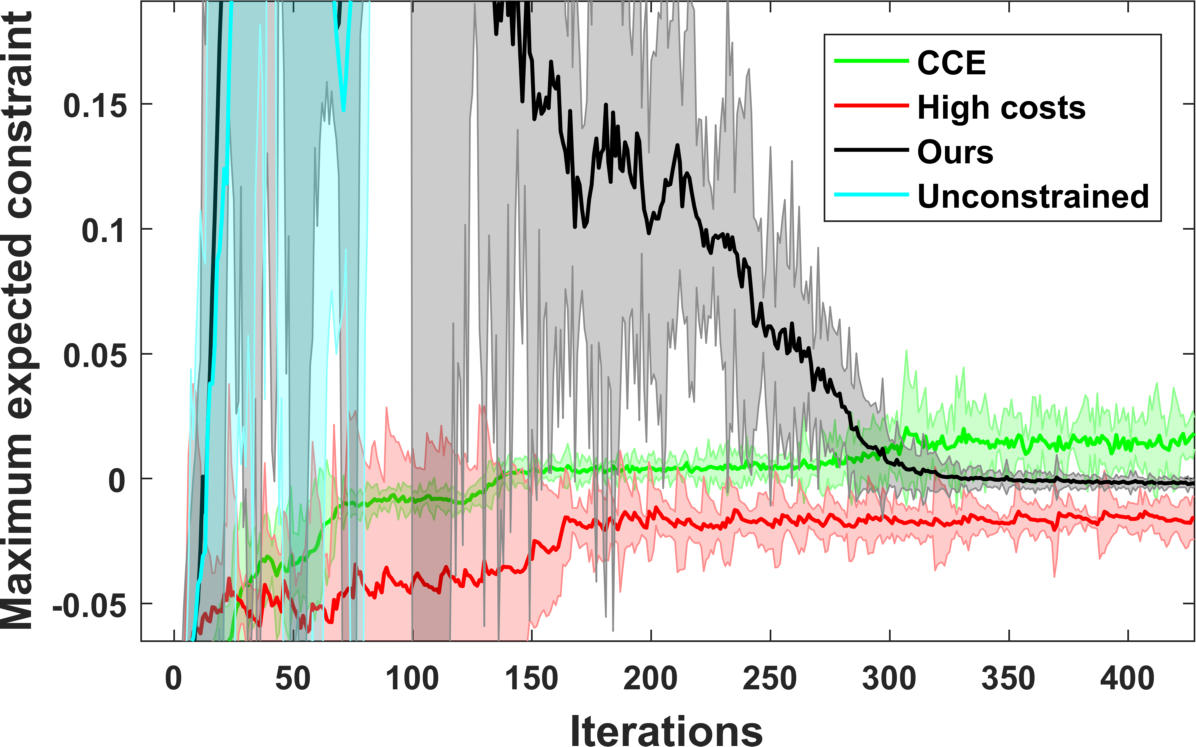}
\caption{Zoomed-in version of Figure 2.}
\end{figure}

\section{Simulations}
\label{section:simulation}

\begin{figure}
\centering
\includegraphics[width=175 pt]{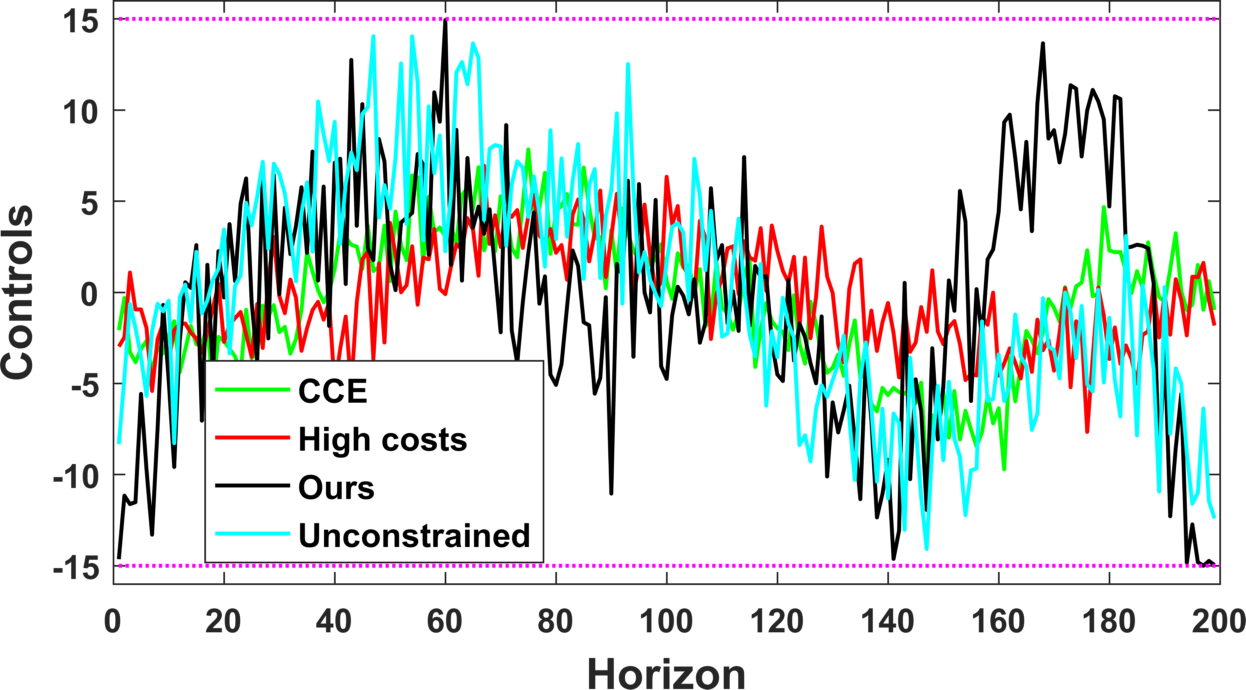}
\caption{Cartpole controls for a single run of the algorithms.}
\end{figure}

\begin{figure}
\centering
\includegraphics[width=198 pt]{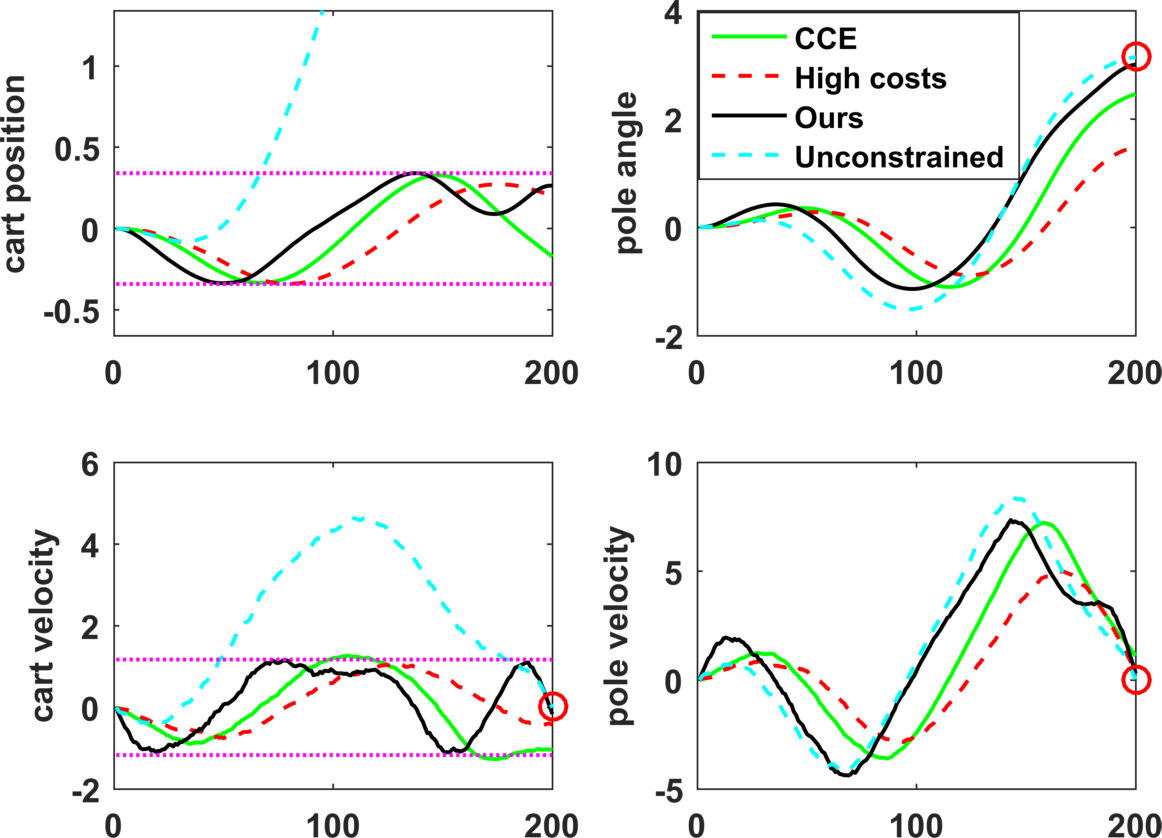}
\caption{Cartpole states for a single run of the algorithms.}
\end{figure}

We compare our proposed scheme against the constrained cross entropy (CCE) method \cite{wen2018constrained} and an algorithm similar to STOMP \cite{kalakrishnan2011stomp} on a cart pole and quadcopter in simulation. The latter algorithm is essentially a path integral control scheme, in which each infeasible rollout will receive a cost equal to $10^4$. This will be denoted by ``High costs'' in our figures. Control constraints for these algorithms were handled via clamping. We will also include the results of the path integral algorithm applied on the uncostrained problem of the cartpole simulation (denoted by ``Unconstrained").

All the simulations are run in Matlab with sample size of 100 rollouts and elite set size of $25$ for CCE. The time steps are 0.01 sec and the task horizon is 200 steps. For the cost plots, we ran all algorithms 4 times and plot the mean and the 97\% confidence region ($\pm3\sigma$ in shaded region).For our algorithm, we will consider the truncated Gaussian distribution for the parameter distribution ${\chi}(\cdot)$. We will also let $\theta^k\equiv u^k$ and only update the mean of the distribution ($\bm{\rho} = \mu \ \text{of}\ {\chi}(\cdot), T(\bm{\theta})=\bm{\theta}$), although the proposed algorithm is applicable to arbitrary distributions in the exponential family and all parameters of the distribution can be updated.

\begin{figure*}[h]
\centering
\includegraphics[width=\linewidth]{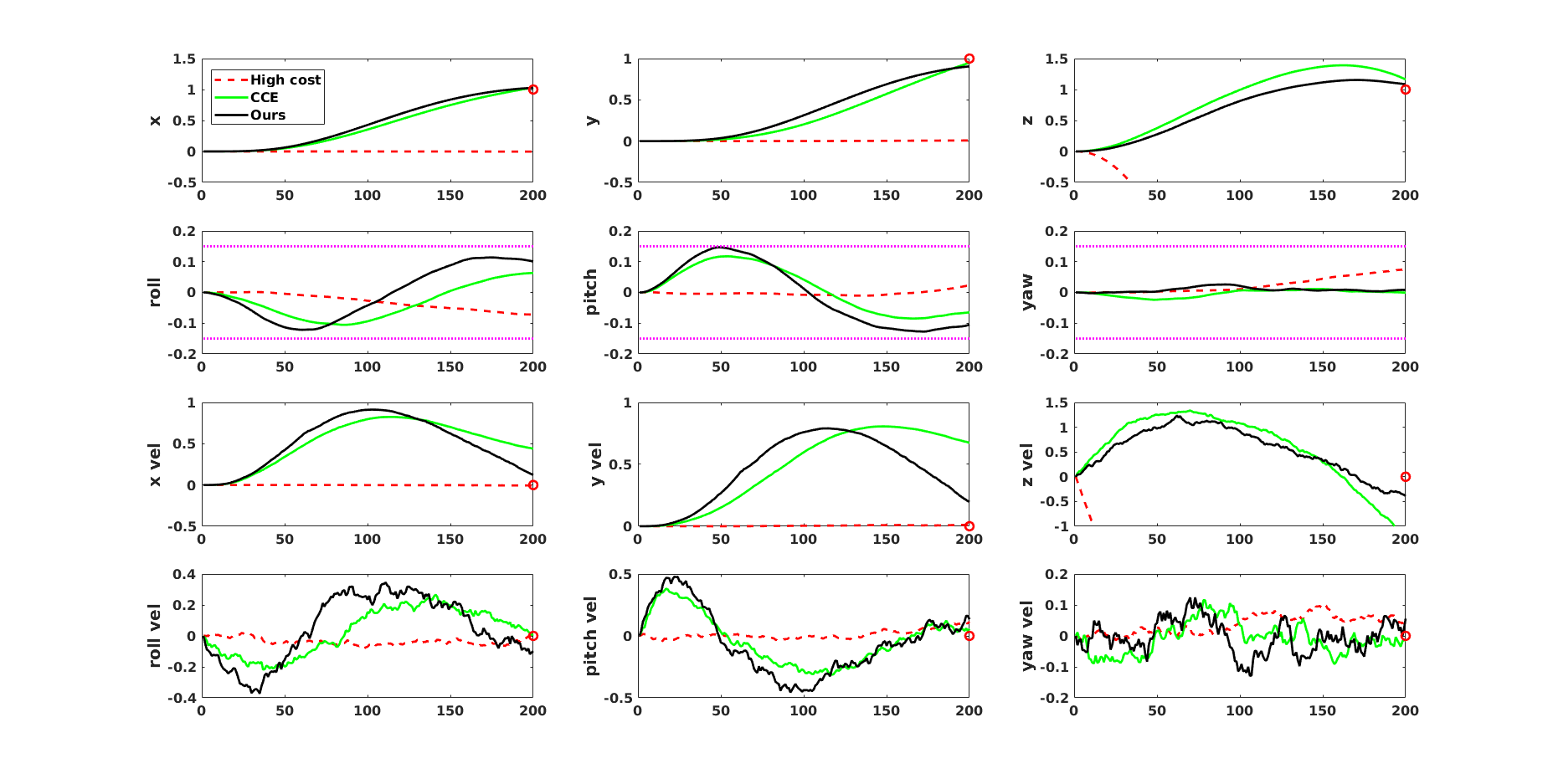}
\caption{Quadcopter states with angular constraints (magenta dotted lines).}
\end{figure*}

\begin{figure}[h]
\centering
\includegraphics[width=200 pt]{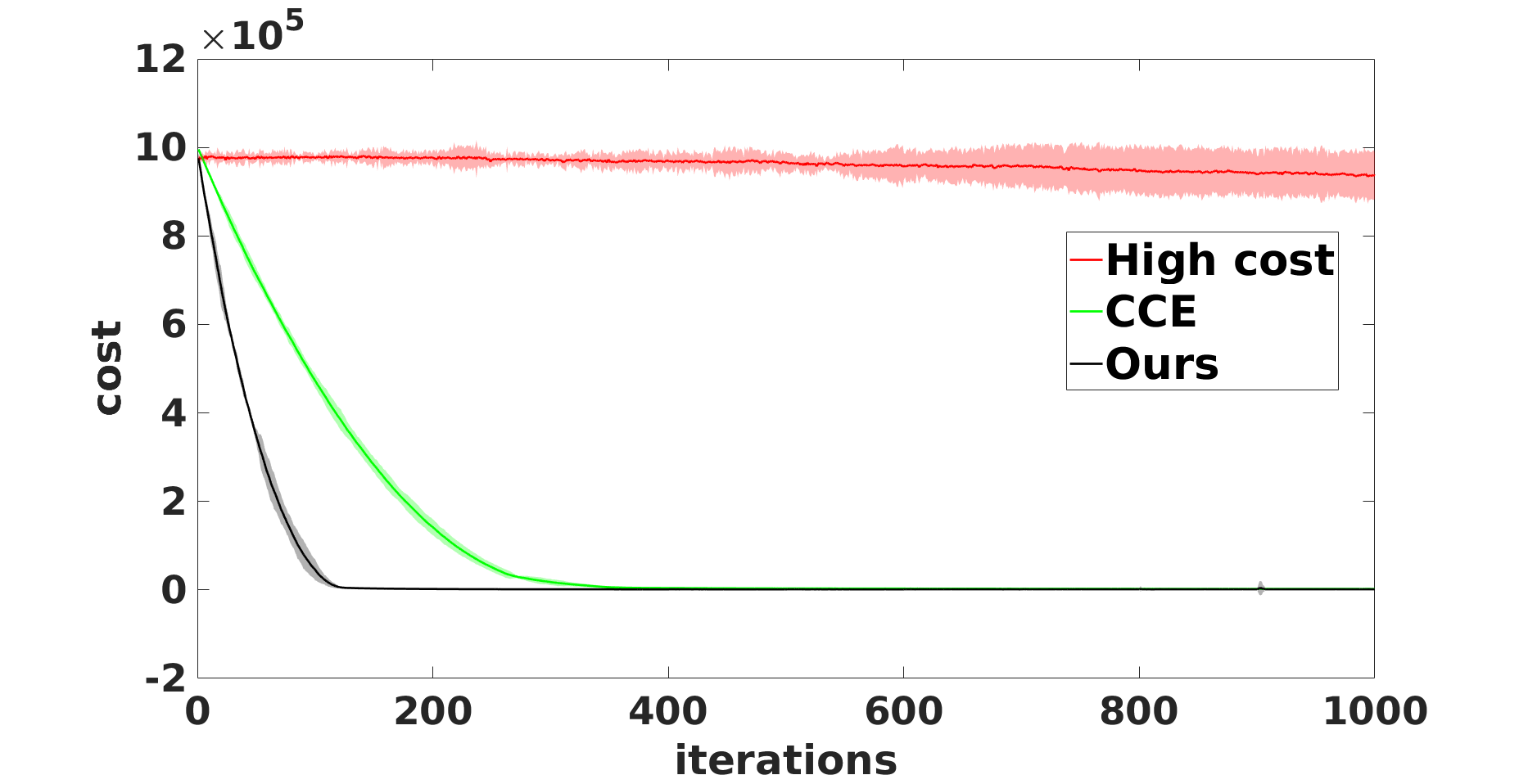}
\caption{Quadcopter cost with a $\pm3\sigma$ confidence interval.}
\end{figure}

\begin{figure}[h]
\centering
\includegraphics[width=\linewidth]{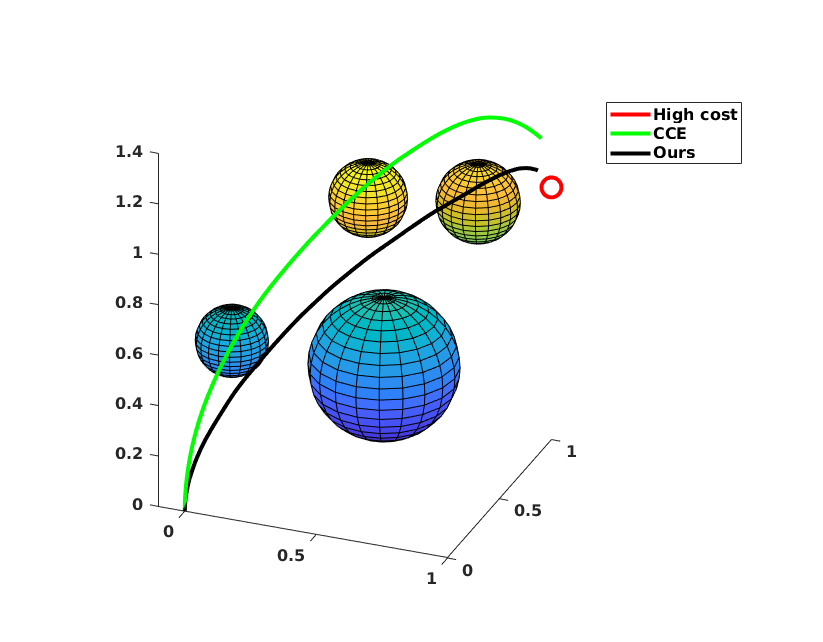}
\caption{Quadcopter state trajectory plot in 3D. The green trajectory denotes the CCE method, and the black trajectory denotes our method. The spheres are the obstacles (nonlinear state constraints).}
\end{figure}

\subsection{Cartpole}
First, we consider the task of a cartpole swing up. The dynamics of the system can be found in \cite{deisenroth2010efficient}. Our initial and target state will be respectively $[0, 0, 0, 0]$, $[-, \pi, 0, 0]$, where the state vector consists of the position of the cart, pole angle, velocity of the cart and pole velocity (note that we do not set a target state for the cart position). The velocity of the pole was influenced by Gaussian noise of variance 0.01. Moreover, the mass of the cart was set to 1 kg, the mass of the pole to 0.1 kg and its length to 0.5 m.

Regarding constraints, we imposed that the position and velocity of the cart satisfy $|\mathsf x^k|\leq 0.34$ $m$ and $|\upsilon^k|\leq 1.17$ $m/s$ for all $k$ respectively. For the controls we considered the box constraints $|u^k|\leq15$. The results of this simulated example can be found in figures 1-5, where we have let all algorithms run for 700 iterations, unless they converged earlier. 

We observe in Fig. 1 that our algorithm and CCE clearly outperform the naive penalization approach (``High costs") within path integral. Moreover, our algorithm outperforms CCE in terms of quality of solutions, while the two methods are comparable in terms of feasibility. In particular, as shown in Figures 2 and 3, our method allows early penalization of the constraints, which eventually gives a better solution. Furthermore, our method typically gives control inputs which are closer to their boundaries, as depicted in Figure 4.

\subsection{Quadcopter}
For the quadcopter, we consider the task of flying from an initial location of $[0,0,0]$ meters to a target location of $[1, 1, 1]$ meters. The dynamics of the system can be found in \cite{michael2010grasp}. We imposed nonlinear state constraints on the quadcopter in the form of four obstacles. Additionally, the quadcopter angles are subject to a maximum angle constraint of 0.15 radians.

Figure 6 compares the states of the quadcopter and demonstrates that our algorithm and CCE again clearly outperform the high cost approach. Additionally, compared to CCE, our method provides a better solution and converges faster. The superior convergence can be observed in figure 7, where our method converges in $1/3$ of the iterations of CCE. Additionally, in figure 8, we can observe that our algorithm provides a more direct trajectory that stays close to the obstacles while satisfying the constraints, whereas the optimal trajectory from CCE avoids the obstacles by going above all of them.



\section{Conclusions}
\label{section:conclusions}
In this paper we proposed a sampling-based algorithm for constrained trajectory optimization problems. Our work was based on the theory of stochastic approximation methods, which allowed us to handle constraints via truncated distributions and exact penalty functions. We presented simulations which demonstrated the effectiveness of our framework over alternative approaches for sampling-based constrained stochastic control. Future directions include having higher-order update formulas, as well as incorporating different methods for handling constraints, such as sequential local approximations with respect to the sampling distributions.








\bibliographystyle{IEEEtran}    
\bibliography{references}

\end{document}